\documentclass{article}

\usepackage{graphicx}

\usepackage{amsmath,amsthm,amssymb}
\numberwithin{equation}{section}
\newtheorem{thm}{Theorem}[section]
\newtheorem{lem}[thm]{Lemma}

\newcommand{\abs}[1]{\left\vert#1\right\vert}

\newcommand{\set}[1]{\left\{#1\right\}}
\newcommand{\Oh}[1]{O\left(#1\right)}

\newcommand{\eps}{\varepsilon}

\newcommand{\To}{\rightarrow}

\newcommand{\inv}{^{-1}}

\newcommand{\var}{\mbox{Var}}

\begin{document}

\title{Optimal Multistage Sampling in a Boundary-Crossing Problem}%

\author{\textbf{Jay Bartroff}\thanks{Address correspondence to Jay Bartroff,
    Department of Statistics, Sequoia Hall, Stanford University,
    Stanford, CA 94305, USA; Fax: 650-725-8977 ; Email:
    bartroff@stat.stanford.edu.  \textbf{Keywords:} Asymptotic;
    Brownian motion; Group sequential; Multistage; Optimality.  \textbf{Subject
    Classifications:} 62L10; 91A20.}\\Department of Statistics, Stanford
    University\\Stanford, CA, USA}%

\date{}

\maketitle

\begin{abstract}Brownian motion with known positive drift is sampled in stages until it crosses a positive boundary $a$.  A family of multistage samplers that control the expected overshoot over the boundary by varying the stage size at each stage is shown to be optimal for large $a$, minimizing a linear combination of overshoot and number of stages.  Applications to hypothesis testing are discussed.\end{abstract}

\section{INTRODUCTION AND SUMMARY}Many problems in theoretical and applied statistics
involve observing a random process until it crosses a
predetermined boundary.  We consider a version of this classical
problem in which Brownian motion $X(t)$, with known drift $\mu>0$ and variance per unit time 1, is sampled in stages until $X(t)\ge a>0$ at the end of a stage.  As an example consider periodic monitoring of a pollutant in a water supply.  There is a critical level for the pollutant above which some action must be taken but below which one will only decide when to test again, basing that decision on the current level.

If one incurs a fixed
cost for each unit sampled and an additional fixed cost for
each stage, then a natural measure of the performance of a
multistage sampler is the sum of these costs upon first crossing the boundary.  In this paper we describe a family of samplers and show they are first-order optimal as $a\To\infty$.

Many aspects of the boundary-crossing or ``first-exit'' problem
are well-studied, though without the multistage aspects considered here.  The powerful methods of renewal theory address
successive exits and the time between such events (see Feller (1971), pages 358-388).  Lorden (1970) obtained sharp, uniform
bounds for the excess over the boundary of random walks.  Siegmund (1985) discusses further applications in sequential analysis.

Schmitz (1993), Cressie (1993), and Morgan (1997) have proved
general existence results for a large class of multistage sampling
problems.  In particular, the theorems of Schmitz show that a optimal sampler does exist for the problem considered
here and that the optimum has the ``renewal-type'' property that
at each stage it behaves as if it were starting from scratch,
given the data so far.  But these authors do not propose specific procedures, and though there is an extensive literature dealing with fully-sequential (one-at-a-time) and group-sequential ($n$-at-a-time) sampling, there have been few investigations of the performance of procedures that vary their sample size from stage to stage.

The families of samplers constructed below, $\delta_{m,h}^o$ and $\delta_{m,z}^+$,
are shown to be first-order asymptotically
optimal in Theorem \ref{thm2.4}.  They have variable stage sizes which
decrease roughly as successive iterations of the function $x\mapsto \sqrt{x\log x}$, while
the average number of stages required is determined by the ratio of the cost per stage to the cost per unit
time in relation to a family of critical functions, $h_m$.  These critical functions define ``critical bands'' -- i.e.,  regions of the first quadrant which are closely related to how close any efficient procedure can be to the boundary after each stage of sampling; Lemma \ref{lem2.5} gives a precise ``in-probability'' lower bound on this distance.  Theorem \ref{thm2.5} then provides a converse statement to the optimality of $\delta_{m,h}^o,
\delta_{m,z}^+$, showing that any competing sampler must use at
least as many stages and follow the same ``schedule'' of $\delta_{m,h}^o$ and $\delta_{m,z}^+$, described in Lemma \ref{lem2.5}.

\section{MULTISTAGE SAMPLERS} Define a \emph{multistage sampling rule} $T$ to be a sequence of nonnegative random variables $(T_1, T_2,\ldots)$ such that, for $k\ge1$ \begin{equation}\label{2.1} T_{k+1}\cdot1\{T_1+\cdots+T_k\le t\}\in \mathcal{E}_t\quad\mbox{for all $t\ge0$,}\end{equation} where $\mathcal{E}_t$ is the class of all random variables determined by $\{X(s): s\le t\}$.  The interpretation of (\ref{2.1}) is that by the time $T^k\equiv T_1+\cdots+T_k$, the end of the first $k$ stages, an observer who knows the values $\{X(s): s\le T^k\}$ also knows the value of $T_{k+1}$, the size of the $(k+1)$st stage.  By a
convenient abuse of notation, we will also let $T$ denote the
total sampling time, $T^M$, where $M\equiv\inf\{m\geq 1:
X(T^m)\geq a\}$, the total number of stages required to cross the boundary $a$.  We will then describe a
\emph{multistage sampler} by the pair $\delta(a)=(T,M)$, where the argument is the initial distance to the boundary.
When there is no confusion as to which sampler is being
used, the shorthand $X_k \equiv X(T^k)$,
$X_0\equiv 0$ will be employed.

Let $c,d>0$ denote the cost per unit time and cost per stage, respectively, and consider the problem of finding the multistage sampler $\delta(a)=(T,M)$ that minimizes $$c\cdot ET+d\cdot EM.$$  Dividing through by $c$, this is seen to be equivalent to minimizing \begin{equation}\label{2.2}ET+h\cdot EM,\end{equation} where $h=d/c$.  By Wald's equation, \begin{equation}\label{2.3}ET=EX(T)/\mu=a/\mu+E(X(T)-a)/\mu\ge a/\mu,\end{equation} so the sampler that minimizes \begin{equation}\label{2.4}E(T-a/\mu)+h\cdot EM\end{equation} also minimizes (\ref{2.2}), and using (\ref{2.4}) instead of (\ref{2.2}) will also lead to a more refined ``first order'' asymptotic theory.

To describe a sampler that asymptotically minimizes
(\ref{2.4}) to first-order, it suffices to consider sequences
$\{(a,h)\}$ such that $a\To\infty$.  We are interested in problems
where optimal procedures use a bounded number of stages and it
turns out that this requires $$h>a^{\eps}$$ for some $\eps>0$.  It
will turn out that good procedures use $m$ stages (almost always)
if, as $a\To\infty$, \begin{equation}\label{2.5} a^{(1/2)^m}(\log
a)^{1/2-(1/2)^m}\ll h\ll a^{(1/2)^{m-1}}(\log
a)^{1/2-(1/2)^{m-1}},\end{equation} where ``$\ll$'' means asymptotically of
smaller order.  We therefore define the \emph{critical functions}
\begin{equation} h_m(x)\equiv x^{(1/2)^m}(\log x)^{1/2-(1/2)^m} \nonumber\end{equation} for
$m=1,2,\ldots$ and $x\ge 1$, with $h_0(x)\equiv x$.  An
essentially complete description of how to achieve asymptotic
optimality is thus given by showing how to proceed in two
cases. The case defined by (\ref{2.5}) is called $\{(a,h)\}$
being in the \emph{$m$th critical band}.  The
other case is
$$h\sim Qh_m(a)$$ for some $Q\in(0,\infty)$, which we refer to as
$\{(a,h)\}$ being on the \emph{boundary between critical bands $m$
and $m+1$}.

It will prove convenient in the sequel to treat $h$ as a function of $a$.  To translate the above formulation into these terms, let $\mathcal{B}_m^o$ be the class of positive functions $h$ such that $\{(a,h(a))\}$ is in the $m$th critical band  (for every sequence of $a$'s approaching $\infty$) and let $\mathcal{B}_m^+$ be the class of positive functions $h$ such that $\{(a,h(a))\}$ is on the boundary between critical bands $m$ and $m+1$  (for every sequence of $a$'s approaching $\infty$).  That is,
\begin{eqnarray*} \mathcal{B}_m^o&\equiv&\set{h:(0,\infty)\To(0,\infty)| \quad h_m(x)\ll h(x)\ll h_{m-1}(x)\quad\mbox{as $x\To\infty$}},\\
\mathcal{B}_m^+&\equiv&\set{h:(0,\infty)\To(0,\infty)| \quad \lim_{x\To\infty}\frac{h(x)}{h_m(x)}\in(0,\infty)},\end{eqnarray*} and let $\mathcal{B}_m=\mathcal{B}_m^o\cup \mathcal{B}_m^+$.  Our notation reflects that, as $a\To\infty$, the average number of stages used by an efficient sampler approaches
$$\begin{array}{cl}
m  &  \mbox{if  $h\in\mathcal{B}_m^o$}  \\
m+\eta  &  \mbox{if  $h\in\mathcal{B}_m^+$},
\end{array}$$ where $\eta\in(0,1)$ is a function of $\lim_{x\To\infty} h(x)/h_m(x)$; Figure 1 summarizes this relationship. Finally, we define the \emph{risk} of a sampler $\delta(a)=(T,M)$ to be \begin{equation} \label{2.6}
R(\delta(a)) = E(T-a/\mu) + h(a)\cdot EM
\end{equation}for a given $h\in\mathcal{B}_m$, some $m\ge1$.  Define the optimal sampler $\delta^*(a)=
(T^*,M^*)$ to be one that achieves $R^*(a) \equiv \inf_{\delta}
R(\delta(a))$. Note that, by (\ref{2.3}), the definition
of risk (\ref{2.6}) is equivalent to the expectation of a linear combination of the so-called ``overshoot,'' $X(T)-a$, and the
number of stages used. 

\begin{table}
\begin{center}
\includegraphics[scale=0.5,angle=0]{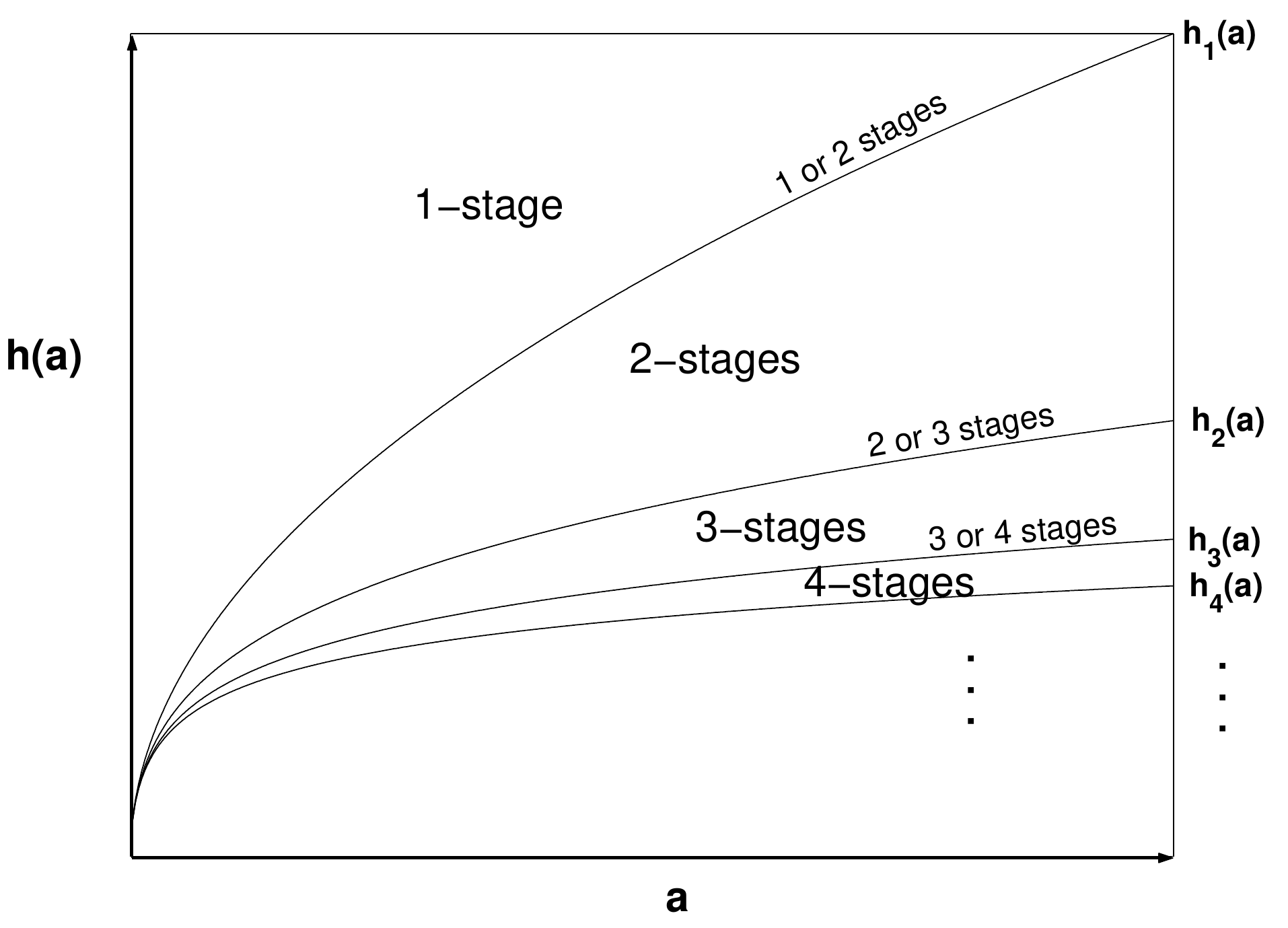}\\
\end{center}
\textbf{\textit{Figure 1.}}  The critical bands and the average number
of stages used by an efficient procedure.
\end{table}

A convenient way of parametrizing stage sizes is by the upper standard normal quantile of the
probability of stopping at the end of the stage.  Thus, for $x>0$ and $z\in\mathbb{R}$ let $t(x,z)$ be the unique solution
of\begin{equation} \label{2.7} \frac{x-\mu
t(x,z)}{\sqrt{t(x,z)}}=z.
\end{equation}  A simple computation gives
\begin{equation}
t(x,z) = x/\mu - \frac{z\sqrt{4x\mu+z^2}-z^2}{2\mu^2}.\nonumber
\end{equation}  Letting $z_p$ denote the upper $p$-quantile of standard normal curve, the probability of being across a boundary $x$ units away at the end of a stage
of size $t(x,z_p)$ is $p$.  An important asymptotic property of $t(x,z)$ is that $$t(x,z)=x/\mu+O(z\sqrt{x})$$ as $x,|z|\To\infty$ if $|z|\ll\sqrt{x}$.

Letting $\Phi$ and $\phi$ denote the standard normal distribution function and density, define $$\Delta(z)\equiv\int_z^{\infty}\Phi(-x)dx=\phi(z)-\Phi(-z)z.$$ The function $\Delta$ will appear often in calculations of expected overshoot or undershoot.  For example, letting $E(Y;A)$ denote $E(Y\cdot 1_A)$, if $Y\sim N(\lambda,\sigma^2)$ then
\begin{eqnarray} E(Y;Y\ge y)&=&\int_y^{\infty}P(Y>w)dw+y\cdot P(Y\ge y)\quad\mbox{integration by parts}\nonumber\\
&=&\int_y^{\infty}\Phi\left(-\frac{w-\lambda}{\sigma}\right)dw+y\cdot \Phi\left(-\frac{y-\lambda}{\sigma}\right)\nonumber\\
&=&\sigma\int_{(y-\lambda)/\sigma}^{\infty}\Phi(-v)dv+y\cdot \Phi\left(-\frac{y-\lambda}{\sigma}\right)\quad\mbox{change of variables}\nonumber\\
&=&\sigma\cdot \Delta\left(\frac{y-\lambda}{\sigma}\right)+y\cdot \Phi\left(-\frac{y-\lambda}{\sigma}\right).\label{2.12}\end{eqnarray} In particular, by taking $y=0$, $\lambda=\mu t(x,z)-x=-z\sqrt{t(x,z)}$, and $\sigma=\sqrt{t(x,z)}$ in (\ref{2.12}) we have
\begin{equation} E[X(t(x,z))-x; X(t(x,z))\geq x] =\sqrt{t(x,z)}\cdot \Delta(z).\label{2.11}\end{equation} Useful asymptotic properties of $\Delta$ are \begin{equation}\label{2.8}\Delta(z)\sim\begin{cases}|z|,&\mbox{as $z\To-\infty$}\\ 
 \phi(z)/z^2,&\mbox{as $z\To+\infty$.}\end{cases}\end{equation} The former is trivial while the latter follows from the classical expansion \begin{equation}\label{2.13}\Phi(-z)=\frac{\phi(z)}{z}\left[1-z^{-2}+O(z^{-4})\right].\end{equation}

\subsection{Geometric Sampling}For $z\in\mathbb{R}$, let $\dot{\delta}_{z}(a)=(T,M)$ be the sampler such that the probability of stopping at the end of each stage is constant at $p\equiv\Phi(-z)$ across the stages, i.e., $$T_k \equiv t(a-X_{k-1},z)\cdot 1\{
X_{k-1}< a\}.$$  Then $M$ is a geometric random variable with mean
$1/p$ and we thus refer to $\dot{\delta}_{z}(a)$ as \emph{geometric sampling}.  Although $p$ is constant across the stages, we
do allow it to vary with $a$, the initial distance to the boundary.

Not only is geometric sampling an interesting random process on its own, but it also has been conjectured that optimal multistage procedures share its stationarity property.  While Theorem \ref{thm2.5} will show this is not true, geometric sampling will prove to be a useful tool for designing the final stages of our optimal samplers in the next section.

Lemma \ref{lem2.1} establishes a fundamental upper bound on the
behavior of geometric sampling and Lemma \ref{lem2.2} gives an
asymptotic bound on the overshoot of $X$ under geometric sampling when
the probability of stopping at each stage approaches 1.  With the
exception of the the main results, Theorems \ref{thm2.4} and
\ref{thm2.5}, the proofs of all theorems and lemmas can be found in
the Appendix.

\begin{lem}\label{lem2.1}  Let $z\in\mathbb{R}$, $q=\Phi(z)$, and \begin{equation}\label{lem2.1.4}g(x) \equiv\frac{\Delta(-z)}{qz}\cdot (x-\mu t(x,z))=
\frac{\Delta(-z)}{2\mu q}(\sqrt{4x\mu+z^2}-z).\end{equation} Then $\dot{\delta}_z(a)=(T,M)$ satisfies
\begin{equation}\label{lem2.1.3}ET-a/\mu\leq \left\{
\begin{array}{lc}
  \frac{q\Delta(z)}{\mu\Delta(-z)}\cdot g(a)+\mu^{-1} \sum _{k\geq 2} g^{(k)}(a)q^{k}, & \mbox{if $z\ge0$} \\
 \frac{q\Delta(z)}{\mu\Delta(-z)}\sum _{k\geq 1}g^{(k)}(a)q^{k-1}, & \mbox{if $z\le0$,} \end{array}\right.\end{equation} where $g^{(k)}$ denotes the $k$th iterate of $g$.\end{lem}

\begin{lem}\label{lem2.2}Let $z(x)$ be a non-positive function such that $\lim_{x\To\infty}z(x)=-\infty$ and $\abs{z(x)}=o(\sqrt{x})$ as $x\To\infty$.  Then $\dot{\delta}_{z(a)}(a)=(T,M)$ satisfies $$EX_M-a=O(\abs{z(a)}\sqrt{a})$$ as $a\To\infty$.\end{lem}

\subsection{The samplers $\delta_{m,h}^o$ and $\delta_{m,z}^+$}
In this section we define two families of samplers that will later be shown to be first order optimal under different conditions.  Namely, the sampler $\delta_{m,h}^o$ will be optimal when $h\in\mathcal{B}_m^o$ and $\delta_{m,z}^+$ will be optimal when $h\in\mathcal{B}_m^+$.

Let $f(x)=(6/\sqrt{\mu})\sqrt{x\log(x+1)}$ and note that $f\inv$ is well-defined since $f$ is increasing.  The family of samplers $\delta_{m,h}^o(a)$ are indexed by a positive integer $m$ and a positive function $h$, and the argument $a$ is the initial distance to the boundary.  We define the family of samplers inductively on $m$ as follows.  Letting $\wedge$ denote min,
\begin{eqnarray*} \delta_{1,h}^o(a)&\equiv& \dot{\delta}_{\zeta(a)}(a),\;\mbox{where $\zeta(x)=-(\sqrt{h(x)}/x^{1/4}\wedge x^{1/7})$.}\\
\delta_{m+1,h}^o(a)&\equiv&\mbox{1st stage $t(a,\sqrt{\log(a/h(a)^2+1)})$,}\\
&&\mbox{followed (if necessary) by $\delta_{m,h\circ f\inv}^o(a-X_1)$.}\end{eqnarray*} 

The family of samplers $\delta_{m,z}^+(a)$ are indexed by a positive integer $m$ and a number $z\in\mathbb{R}$.  They are defined inductively on $m$ as follows:   
\begin{eqnarray*} \delta_{1,z}^+(a)&\equiv&\mbox{1st stage $t(a,z)$, followed (if necessary) by}\\
&&\mbox{$\dot{\delta}_{\nu(a-X_1)}(a-X_1)$, where $\nu(x)=-\sqrt{\log(x+1)}$.}\\
\delta_{m+1,z}^+(a)&\equiv& \mbox{1st stage $t(a,\sqrt{(1-2^{-m})\log(a+1)})$,}\\
&&\mbox{followed (if necessary) by $\delta_{m,z}^+(a-X_1)$.}\end{eqnarray*}

The main theorems in this section establish the operating characteristics of $\delta_{m,h}^o$ and $\delta_{m,z}^+$.  Theorem \ref{thm2.1} will show that $\delta_{m,h}^o$ uses $m$ stages almost always and has overshoot bounded by $h$, asymptotically.  Theorem \ref{thm2.2} will show that $\delta_{m,z}^+$ uses $m$ or $m+1$ stages almost always, and that the probability of using $m$ stages has upper quantile approaching $z$.

\begin{thm}\label{thm2.1}If $h\in\mathcal{B}_m^o$ then $\delta_{m,h}^o(a)=(T,M)$ satisfies
\begin{eqnarray} ET-a/\mu&=&o(h(a))\label{thm2.1.1}\\
EM&\To&m\label{thm2.1.2}\end{eqnarray} as $a\To\infty$.\end{thm}

Before establishing the operating characteristics of $\delta_{m,z}^+$ in Theorem \ref{lem2.2}, we introduce a family of positive constants that appear in the first order overshoot of $\delta_{m,z}^+$.  For $m\ge 1$ define \begin{equation}\label{2.9}\kappa_m=\kappa_m(\mu) = \mu^{-2+(1/2)^m}
\prod_{i=1}^{m-1} [(1/2)^{m-1-i}-(1/2)^{m-1}]^{(1/2)^{i+1}}.\end{equation} The next lemma contains the key property of the $\kappa_m$ that we will need.

\begin{lem}\label{lem2.3}For $m\ge 1$, as $x\To\infty$ $$\kappa_mh_m(\sqrt{(1-2^{-m})/\mu\cdot x\log x})\sim\kappa_{m+1}h_{m+1}(x).$$\end{lem}

We adopt the notation $F\lesssim G$ to denote $F\le (1+o(1))\cdot G$.  Theorem \ref{thm2.2} establishes the operating characteristics of $\delta_{m,z}^+$.

\begin{thm}\label{thm2.2}Let $m\in\mathbb{N}$ and $z\in\mathbb{R}$.  Then $\delta_{m,z}^+(a)=(T,M)$ satisfies
\begin{eqnarray} ET-a/\mu&\lesssim& \Delta(z)\kappa_m h_m(a)\label{thm2.2.1}\\
EM&\To&m+\Phi(z)\label{thm2.2.2}\end{eqnarray} as $a\To\infty$.\end{thm}

\subsection{Optimality of $\delta_{m,h}^o$ and $\delta_{m,z}^+$}\label{sec2.1.3}

In this section we state our main optimality results of the paper, Theorems \ref{thm2.4} and \ref{thm2.5}.  Theorem \ref{thm2.4} shows that $\delta_{m,h}^o$ is first-order optimal when $h\in\mathcal{B}_m^o$ and that $\delta_{m,z^*}^+$ is first-order optimal when $h\in\mathcal{B}_m^+$, where the value $z^*$ is determined by $\lim_{x\To\infty} h(x)/h_m(x)$.  Theorem \ref{thm2.5} provides a converse to Theorem \ref{thm2.4}, showing that any efficient sampler must behave like $\delta_{m,h}^o$ and $\delta_{m,z^*}^+$.

\begin{thm}\label{thm2.4}If $h\in\mathcal{B}_m^o$, then \begin{equation}
\label{thm2.4a}
R(\delta_{m,h}^o(a))\sim m\cdot h(a)\sim R^*(a)
\end{equation}as $a\To\infty$.  If $h\in\mathcal{B}_m^+$, then
\begin{equation}
\label{thm2.4f}
R(\delta_{m,z^*}^+(a))\sim \left[m+\Phi(z^*)+\Delta(z^*)\cdot \frac{\phi(z^*)}{1-\Phi(z^*)}\right] h(a)\sim R^*(a)
\end{equation} as $a\To\infty$, where $z^*$ is the unique solution of \begin{equation}\label{thm2.4.1}\frac{\phi(z^*)}{1-\Phi(z^*)} = \lim_{x\To\infty}\frac{\kappa_m h_m(x)}{h(x)}.\end{equation}
\end{thm}

Before proving Theorem \ref{thm2.4} we introduce a family of iterated functions $F_h^{(k)}$ which are the order of magnitude of the best-possible undershoot after $k$ stages, $a-X_k$.  Lemma \ref{lem2.5} makes this explicit and Lemma \ref{lem2.4} establishes a link between the $F_h^{(k)}$ and the critical functions $h_k$.

For $x\geq y^2>0$ define $F_y(x) = \sqrt{x\log (x/y^2)}$.  For a positive function $h$ and $x$ such that $h^2(x)\leq x$
define $$F_{h}^{(k)}(x) \equiv F_y^{(k)}(x)|_{y=h(x)}.$$  Note that
$h(\cdot)$ is not iterated, e.g., 
\begin{eqnarray*} F_h^{(2)}(x)&=&\left[x\log(x/h^2(x))\right]^{1/4}\sqrt{\log\left(\frac{\sqrt{x\log(x/h^2(x))}}{h^2(x)}\right)}\\
&\ne& \left[x\log\left(x/h^2(\sqrt{x\log(x/h^2(x))})\right)\right]^{1/4} \sqrt{\log\left(\frac{\sqrt{x\log(x/h^2(x))}}{h^2(\sqrt{x\log(x/h^2(x))})}\right)}.\end{eqnarray*}

Lemma \ref{lem2.5} establishes an ``in probability'' lower bound of the order $F_h^{(k)}$ on how close to the boundary any efficient sampler can be after each of the first $m-1$ stages when $h\in\mathcal{B}_m$.

\begin{lem}\label{lem2.5} If $h\in\mathcal{B}_m$ and $\delta$ is any sampler such that $R(\delta(a))=O(h(a))$, then for any $\eps>0$ and $0\le k\le m-1$
\begin{equation}P\left(a-X_k\ge (1-\eps)\cdot (1/\mu)^{1-2^{-k}}F_{h}^{(k)}(a)\right)\To 1\label{lem2.5.0}
\end{equation} as $a\To\infty$.\end{lem}

Lemma \ref{lem2.4} shows that, when $h\in\mathcal{B}_m$, square roots of the iterates $F_{h}^{(k-1)}(a)$ are roughly constant multiples of the critical functions $h_k$.  The constants themselves are given by the solutions of the following recurrence relation.  For $1\le k\le m$ define $C_k^m$ to be the unique solution of \begin{equation}\label{lem2.4.6}C_{k+1}^m = \sqrt{C_k^m}\cdot[(1/2)^{k-1}-(1/2)^{m-1}]^{1/4};\quad C_1^m=1.\end{equation} After taking logarithms, solving (\ref{lem2.4.6}) amounts to solving a difference equation.  This computation gives \begin{equation}\label{lem2.4.3}C_k^m=\prod_{i=1}^{k-1}\left[(1/2)^{k-1-i}-(1/2)^{m-1}\right]^{(1/2)^{i+1}},\end{equation} where it is understood that an empty product equals 1.  Note also that \begin{equation}\label{2.9}\kappa_m = (1/\mu)^{2-(1/2)^m}C_m^m.\end{equation}

\begin{lem}\label{lem2.4} If $h\in\mathcal{B}_m^+$, then \begin{equation}
\label{lem2.4.2}\sqrt{F_{h}^{(k-1)}(x)}\sim C_k^mh_k(x)\quad\mbox{as $x\To\infty$, for $1\leq k\leq m$.}
\end{equation} If $h\in\mathcal{B}_m^o$, then
\begin{equation} \label{lem2.4.1} C_k^{m-1}\lesssim\frac{\sqrt{F_{h}^{(k-1)}(a)}}{h_k(x)}\lesssim C_k^m \quad\mbox{as $x\To\infty$, for $1\leq k<m$.}
\end{equation}\end{lem}

\noindent {\bf Proof of Theorem \ref{thm2.4}.}  Assume that $h\in\mathcal{B}_m^o$.  The left hand side of (\ref{thm2.4a}) holds by Theorem \ref{thm2.1}. Now
\begin{eqnarray*} R^*(a)&\le&R(\delta_{m,h}^o(a))\quad\mbox{the Bayes property}\\
&=&O(h(a))\end{eqnarray*} by Theorem \ref{thm2.1}, so Lemma \ref{lem2.5} applies to $\delta^*(a)$.  Letting $X^*_k$ denote the $\delta^*$-sampled process,
\begin{eqnarray*} R^*(a) &\geq&  h(a)EM^*\\
&\geq& h(a)mP(M^*\geq m)\\
&=&h(a)mP(a-X^*_{m-1}> 0)\\
&\geq& h(a)m P\left(a-X^*_{m-1}\geq (1/2)(1/\mu)^{1-2^{-(m-1)}}
F_{h}^{(m-1)}(a)\right)\\
&\sim& h(a)m\cdot 1\end{eqnarray*} by Lemma \ref{lem2.5}, proving (\ref{thm2.4a}).

If $h\in\mathcal{B}_m^+$ then  $\lim_{x\To\infty}h(x)/h_m(x)$ is positive and finite, hence $$\lim_{x\To\infty}\frac{\kappa_mh_m(x)}{h(x)}=\kappa_m\left[\lim_{x\To\infty} \frac{h(x)}{h_m(x)}\right]\inv$$ is positive and finite as well.  The function $$z\mapsto\frac{\phi(z)}{1-\Phi(z)}$$ increases from 0 to $\infty$ as $z$ ranges from $-\infty$ to $\infty$, so the equation (\ref{thm2.4.1}) has a unique solution, $z^*$. Theorem \ref{thm2.2} shows that
\begin{eqnarray*} R(\delta_{m,z^*}^+(a))&\lesssim&\Delta(z^*)\kappa_mh_m(a)+h(a)[m+\Phi(z^*)]\\
&\sim&\Delta(z^*)\frac{\phi(z^*)}{1-\Phi(z^*)}h(a)+h(a)[m+\Phi(z^*)]\quad\mbox{by definition of $z^*$}\\
&=&\left[m+\Phi(z^*)+\Delta(z^*)\frac{\phi(z^*)}{1-\Phi(z^*)}\right]h(a).\end{eqnarray*} Suppose $\eps>0$.  Again $R^*(a)\leq R(\delta_{m,z^*}^+(a))=O(h(a))$, so by Lemma \ref{lem2.5} \begin{equation}\label{thm2.4.3}P(a-X^*_{m-1}\geq (1-\varepsilon)(1/\mu)^{1-2^{-(m-1)}} F_{h}^{(m-1)}(a))\To 1.\end{equation} Let $(T^{*(m)}, M^{*(m)})$ denote the continuation of $\delta^*$ after the $(m-1)$st stage, i.e.,
\begin{eqnarray*}
M^{*(m)} & = & [M^*-(m-1)]^+, \\
T^{*(m)} & = & T^* - (T_1^*+\cdots + T_{m-1}^*),
\end{eqnarray*}  and for $y>0$ define $$\varphi(y) = E[\mu^{-1}(X_{M^{*(m)}}-y) + h(a)M^{*(m)}| a-X^*_{m-1} = y].$$  We will show below that $\varphi(y)$ is non-decreasing in $y$.  Let $$\gamma = (1-\varepsilon)(1/\mu)^{1-2^{-(m-1)}} F_{h}^{(m-1)}(a).$$  We now compute a lower bound for $\varphi(\gamma)$.  Letting $$p=P(M^{*(m)}=1|a-X^*_{m-1} = \gamma),$$ we have
\begin{eqnarray}
\lefteqn{\mu^{-1}E(X_{M^{*(m)}}-(a-X^*_{m-1})| a-X^*_{m-1} = \gamma)}\nonumber\\
&&\geq\mu^{-1}E[(X_{M^{*(m)}}-(a-X^*_{m-1}))1\{M^{*(m)}=1\}| a-X^*_{m-1} = \gamma]\nonumber\\
&&= \mu^{-1} \Delta(z_p)\sqrt{t(\gamma,z_p)}\nonumber\\
&&\sim \mu^{-1}\Delta(z_p)\sqrt{\gamma/\mu}\nonumber\\
&&=\mu\inv \Delta(z_p)\sqrt{(1-\eps)(1/\mu)^{2-2^{-m+1}}F^{(m-1)}_{h}(a)}\nonumber\\
&&\sim\Delta(z_p)\sqrt{1-\eps}\cdot (1/\mu)^{2-2^{-m}}C_m^mh_m(a)\quad\mbox{by Lemma \ref{lem2.4}}\nonumber\\
&&=\Delta(z_p)\sqrt{1-\eps}\cdot \kappa_mh_m(a)\nonumber\quad\mbox{by (\ref{2.9})} \\
&&\sim\Delta(z_p)\sqrt{1-\eps}\cdot\frac{\phi(z^*)}{1-\Phi(z^*)}h(a),\label{thm2.4b}
\end{eqnarray}by definition of $z^*$.  Also, $E(M^{*(m)}|a-X^*_{m-1}=\gamma)\geq 2-p$, and combining this with (\ref{thm2.4b}) gives, for sufficiently large $a$, 
\begin{eqnarray} \varphi(\gamma)&\geq& \left[\Delta(z_p)\sqrt{1-\eps}\cdot\frac{\phi(z^*)}{1-\Phi(z^*)} + 2-p\right]h(a) \cdot(1-\eps)\nonumber\\
&\ge&\left[\Delta(z_p)\frac{\phi(z^*)}{1-\Phi(z^*)} + 2-p\right]h(a) \cdot(1-2\eps).\label{thm2.4.2}\end{eqnarray} Letting $Y = a-X^*_{m-1}$ and $V = \{Y\geq \gamma\}$, we have
\begin{eqnarray}
R^*(a) & = & E[\mu^{-1}(X_{M^*}-a)+ h(a)M^*]\nonumber \\
&\geq&E[\mu^{-1}(X_{M^{*(m)}}-Y)+ h(a)(m-1+M^{*(m)});V]\nonumber \\
&\geq&E[\varphi(Y)+ (m-1)h(a);V]\nonumber\\
&\geq& [\varphi(\gamma)+ (m-1)h(a)]P(V)\quad\mbox{($\varphi$ non-decreasing)}\nonumber\\
&\ge&\left[\Delta(z_p)\frac{\phi(z^*)}{1-\Phi(z^*)} + m+1-p\right]h(a) \cdot(1-2\eps)\cdot P(V)\quad\mbox{by (\ref{thm2.4.2})}\nonumber\\
&\ge& \left[\Delta(z_p)\frac{\phi(z^*)}{1-\Phi(z^*)} + m+1-p\right]h(a) \cdot(1-3\eps)\label{thm2.4e}
\end{eqnarray} for sufficiently large $a$ since $P(V)\To1$ by (\ref{thm2.4.3}).  Using basic calculus, it can be shown that the expression in brackets in (\ref{thm2.4e}) achieves its unique minimum when $p=1-\Phi(z^*)$, hence $$R^*(a)\ge \left[\Delta(z^*)\frac{\phi(z^*)}{1-\Phi(z^*)} + m+\Phi(z^*)\right]h(a) \cdot(1-3\eps)$$ for sufficiently large $a$.  Since $\eps$ was arbitrary, this completes the proof of (\ref{thm2.4f}) and hence the theorem once we verify that $\varphi(\cdot)$ is non-decreasing.  

Fix $a>0$ and let $0<y\leq y'$.  Let $(T^{\prime (m)},M^{\prime (m)})$ denote the continuation of $\delta^*$ after the $(m-1)$st stage that uses the same probability of being over the boundary at the end of each stage as $(T^{*(m)}, M^{*(m)})$ when starting from $a-X^*_{m-1}=y'$.  Then \begin{equation}
\label{thm2.4g}
E(M^{\prime (m)}|a-X^*_{m-1}=y) = E(M^{*(m)}|a-X^*_{m-1}=y')
\end{equation}and, letting $$p_1 = P(M^{*(m)}=1|a-X^*_{m-1}=y')= P(M^{\prime (m)}=1|a-X^*_{m-1}=y),$$ we have
\begin{equation*}\begin{split}
E[(X_{M^{\prime (m)}}-y)1\{M^{\prime (m)}=1\}| &a-X^*_{m-1}=y]= \Delta(z_{p_1})\sqrt{t(y,z_{p_1})}\\
&\leq  \Delta(z_{p_1})\sqrt{t(y',z_{p_1})}\quad\mbox{since $y\leq y'$}\\
&=E[(X_{M^{*(m)}}-y')1\{M^{*(m)}=1\}| a-X^*_{m-1}=y'].\end{split}\end{equation*}  Similar
arguments inductively give
$$E[(X_{M^{\prime (m)}}-y)1\{M^{\prime (m)}>1\}| a-X^*_{m-1}=y]\leq E[(X_{M^{*(m)}}-y')1\{M^{*(m)}>1\}| a-X^*_{m-1}=y'],$$ and these last two bounds show \begin{equation}
E(X_{M^{\prime (m)}}-y| a-X^*_{m-1}=y)\leq E(X_{M^{*(m)}}-y'| a-X^*_{m-1}=y').\label{thm2.4c}
\end{equation}Then \begin{eqnarray*}
\varphi(y)&\leq& E[\mu^{-1}(X_{M^{\prime (m)}}-y)+h(a)M^{\prime (m)}|a-X^*_{m-1}=y]\quad\mbox{(optimality of $(T^{*(m)}, M^{*(m)})$)}\\
&\leq& E[\mu^{-1}(X_{M^{*(m)}}-y')+h(a)M^{*(m)}|a-X^*_{m-1}=y']\quad\mbox{(by (\ref{thm2.4g}) and (\ref{thm2.4c}))}\\
&=&\varphi(y'),
\end{eqnarray*}finishing the proof.\qed\bigskip

 The final result of this section is a type of converse to Theorem \ref{thm2.4}, showing that good samplers must behave like $\delta_{m,h}^o, \delta_{m,z}^+$ in not only the sense that $m$ stages are necessary when $h\in\mathcal{B}_m$, but also that any efficient sampler must follow the same ``schedule'' that $\delta_{m,h}^o$ and $\delta_{m,z}^+$ follow for the first $m-1$ stages, described in Lemma \ref{lem2.5}.

 \begin{thm}\label{thm2.5}Assume that $h\in\mathcal{B}_m$ and let $$
\delta(a) = \left\{ \begin{array}{ll}
  \delta_{m,h}^o(a)    &  \mbox{if $h\in\mathcal{B}_m^o$}  \\
  \delta_{m,z^*}^+(a)    &   \mbox{if $h\in\mathcal{B}_m^+$,}
\end{array}\right.$$ where $z^*$ is as in (\ref{thm2.4.1}). If $\delta'(a)=(T,M)$ is a sampler such that there is a sequence $a_i\To\infty$ with
\begin{equation}
\label{thm2.5b} P\left(a_i-X_k\ge
(1-\varepsilon)(1/\mu)^{1-2^{-k}}F_{h}^{(k)}(a_i)\right)
\quad\mbox{bounded below 1}
\end{equation}for some $1\leq k<m$ and $\varepsilon>0$, then \begin{equation}
\label{thm2.5a} \frac{R(\delta'(a_i))-R(\delta(a_i))}{h(a_i)}\To+\infty
\end{equation} as $i\To\infty$.  In particular, (\ref{thm2.5a}) holds if $P(M\geq m)\not\To 1$.\end{thm}

\bigskip\noindent{\bf Proof.}  If (\ref{thm2.5b}) holds then Lemma \ref{lem2.5} implies that $$\frac{R(\delta'(a_i))}{h(a_i)}\To\infty$$ as $i\To\infty$.  We know that $R(\delta(a_i))=O(h(a_i))$ by Theorem \ref{thm2.4}, so $$\frac{R(\delta'(a_i))-R(\delta(a_i))}{h(a_i)}=\frac{R(\delta'(a_i))}{h(a_i)}-O(1)\To+\infty$$ as $i\To\infty$, proving the first claim. 

If $$P\left(a-X_k\ge
(1-\varepsilon)(1/\mu)^{1-2^{-k}}F_{h}^{(k)}(a)\right)\To1$$ for all $1\leq k<m$, then \begin{eqnarray*} P(M\geq m)&=&P(a-X_{m-1}>0)\\
&\ge& P\left(a-X_{m-1}\ge
(1-\varepsilon)(1/\mu)^{1-2^{-(m-1)}}F_{h}^{(m-1)}(a)\right)\\
&\To& 1,\end{eqnarray*} proving the second assertion.\qed

\section{AN APPLICATION TO HYPOTHESIS TESTING}\label{sec5}

 In this section we discuss how the above multistage samplers can be used to construct efficient multistage hypothesis tests and give a numerical example of the performance of such a test.  Let $Y_1,Y_2,\ldots$ be i.i.d.\ with density $f$ and consider testing the simple hypotheses

\begin{equation}\label{app2}H_0: f=f_0\quad\mbox{vs.}\quad H_1: f=f_1\end{equation} in stages.  We can describe multistage tests of these hypotheses by triples $(N,M,D)$, where $N$ is the

total number of observations, $M$ is the total number of stages

used, and $D$ is the decision variable, taking values in

$\{0,1\}$.  One measure of the performance of $(N,M,D)$ in testing $H_0, H_1$ is the \textit{integrated risk}, which we define as

\begin{equation} \label{app1} r(N,M,D)=\sum_{i=0}^1 [cE_iN
  +dE_iM+w_iP_i(D\ne i)]\pi_i,
\end{equation}where $0<c,d<1$ represent the cost per observation and cost per stage,

$\pi_i$ is the prior distribution on $\{f_0, f_1\}$, and $w_i>0$

represent the penalty for a wrong decision.  Here $E_i$ and $P_i$ denote expectation and probability under $f_i$.

The multistage samplers discussed above, with some simple modifications (like making sure each stage size is an integer), tell us how to sample the $Y_j$ in stages by observing some random process until it crosses a boundary.  The relevant random process here, taking the place of the Brownian motion above, is the log-likelihood process, which we now define.  

Assume that $$E_i
\log^2\left(\frac{f_i(Y_1)}{f_{1-i}(Y_1)}\right)<\infty\quad\mbox{for
  $i=0,1$}$$ and let $$\sigma_i^2 = \var_i\log\left(\frac{f_i(Y_1)}{f_{1-i}(Y_1)}\right),\quad \mu_i =\sigma_i\inv E_i\log\left(\frac{f_i(Y_1)}{f_{1-i}(Y_1)}\right).$$ Define the log-likelihood process

$$X_i(n) =  \sigma_i\inv \sum_{j=1}^n
\log\left(\frac{f_i(Y_j)}{f_{1-i}(Y_j)}\right).$$ Note that, like the
Brownian motion above, $E_iX_i(n)=n\mu_i>0$ and $\var_iX_i(n)=n$.  Thus, we can use a given multistage sampler $\delta'$ with boundary $a$ to sample the $Y_j$ by treating $f_i$ as the true, underlying density and sampling according to $\delta'$ until $X_i(n)\ge a$ at the end of a stage.  If $a$ is large, this is compelling evidence that $f_i$ is indeed the true density, and thus provides us with a ``one decision'' test of the hypotheses (\ref{app2}).  We can construct an ordinary ``two decision'' test from two samplers, $\delta_0$ and $\delta_1$, by somehow choosing a first stage size, computing the maximum likelihood estimated of the true hypothesis $\hat{\imath}\in\{0,1\}$ from the data observed in that first stage, and then continuing with $\delta_{\hat{\imath}}$ and a slightly modified stopping rule (e.g., ``stop when $\abs{X_i(n)}\ge a$'') to protect against an error in $\hat{\imath}$.  A natural choice for the size of the first stage is the minimum of the two stage sizes dictated by $\delta_0,\delta_1$.

Let $\delta$ be the test constructed in this manner from the two samplers $\delta_{m_i^*, d/c}^o(a_i)$, $i=0,1$, where $a_i\equiv\sigma_i\inv\log d\inv$ and $$m_i^*\equiv\inf \set{m\ge 1:\kappa_m(\mu_i)h_m(a_i)\le d/c\le \kappa_{m+1}(\mu_i)h_{m+1}(a_i)}.$$  Bartroff (2004, 2005) shows that, under general conditions on the $f_i$, $\delta$ minimizes the integrated risk (\ref{app1}) to second order as $c,d\To0$ at specified rates.  Bartroff (2004) and a third paper by the author in this series will extend these results to composite hypotheses about the parameter of an exponential family.  See these references for more details on the testing problem.

\begin{table}

\noindent \textbf{Table 1.}  Results for testing $\mu = .25$ vs. $\mu = -.25$
about a Gaussian mean, with $d = .001$, $\pi_i = 1/2$, $w_i = 1$.

\begin{center}

\begin{tabular}{c c r c c}



\hline

Procedure & $EN$ & $EM$ & $r$ & $r(\delta)/r$ (\%)\\ \hline

\multicolumn{5}{c}{$d/c = 1$}\\

$\delta$ & 62.2 & 5.2 & .068 & 100\\

$\delta_g(1)$ & 57.5 & 57.5 &.115 & 59.1\\

$\delta_g(15)$ & 64.9 & 4.6 &.073 & 92.9\\

$\delta_g(30)$ & 76.7 & 2.6 & .080& 85.0\\

\hline

\multicolumn{5}{c}{$d/c = 5$}\\

$\delta$ & 68.3 & 2.9 & .017 & 100\\

$\delta_g(1)$ & 57.5 & 57.5 &.070 & 24.0\\

$\delta_g(22)$ & 72.7 & 3.3 &.018 & 92.8\\

$\delta_g(44)$ & 83.6 & 1.9 & .019& 88.4\\

\hline

\multicolumn{5}{c}{$d/c = 10$}\\

$\delta$ & 76.6 & 1.9 & .0097 & 100\\

$\delta_g(1)$ & 57.5 & 57.5 &.0644 & 15.1\\

$\delta_g(37)$ & 80.5 & 2.2& .0104 & 93.3\\

$\delta_g(74)$ & 97.6 & 1.3 & .0112& 86.6\\ \hline

\end{tabular}

\end{center}
\end{table}

Table 1 contains the results of a numerical experiment comparing
$\delta$ with group-sequential (i.e., constant stage-size) testing of
the hypotheses $\mu=.25$ vs. $\mu=-.25$, about the mean of Gaussian
random variables with unit variance.  $\delta_g(k)$ denotes
group-sequential testing with constant stage-size $k$, which samples
until $$\abs{\sum_j \log\left(\frac{f_1(Y_j)}{f_0(Y_j)}\right)}\geq
\log d\inv$$ at the end of a stage, which is equivalent to the
stopping rule of $\delta$.  For each value of $d/c$, the operating
characteristics of $\delta_g(k)$ are given for three values of $k$:
$k=1$, the best possible $k$ (determined by simulation), and two times
the best possible $k$. The results show significant improvement in the
integrated risk of the variable stage-size test $\delta$ upon
$\delta_g$. This improvement decreases for large values of $d/c$, but this is to be expected since the number of stages of any reasonable test will approach 1 in this limit.

Here we constructed tests from the samplers $\delta_{m,h}^o$.  In practice, tests constructed from the samplers $\delta_{m,z}^+$ also perform well and behave almost identically as those constructed from $\delta_{m,h}^o$.  The choice here was made merely to simplify presentation.

\appendix
\section{PROOFS OF THEOREMS AND LEMMAS}
\noindent {\bf Proof of Lemma \ref{lem2.1}.} First we will prove \begin{equation}\label{lem2.1.1}E(a-X_k | M>k)\leq
g^{(k)}(a)\quad\mbox{ for all $k\geq 0$.}\end{equation}  The $k=0$ case is trivial and we have
\begin{eqnarray*} E(a-X_{k+1}| M> k+1, X_k) &=&E[(a-X_k)-(X_{k+1}-X_k)|M> k+1, X_k]\\
&=&\Delta(-z)\sqrt{t(a-X_k,z)}/q\\
&=&\Delta(-z)\frac{(a-X_k)-\mu t(a-X_k,z)}{qz}\\
&=& g(a-X_k).\end{eqnarray*}
$g$ is increasing and concave, so by Jensen's inequality and the
induction hypothesis
\begin{eqnarray} E(a-X_{k+1}|M>k+1)&=&E(g(a-X_k)|M>k+1)\nonumber\\
&\leq& g(E(a-X_k|M>k+1))\nonumber\\
&=& g(E(a-X_k|M>k))\label{lem2.1.2}\\
&\leq& g(g^{(k)}(a))=g^{(k+1)}(a),\nonumber\end{eqnarray} proving (\ref{lem2.1.1}).  In (\ref{lem2.1.2}) we use that $E(a-X_k|M>k+1) = E(a-X_k|M>k)$; this is
true since the value of $X_k$ and the number of additional stages
required to cross the boundary are independent, as long as $X_k <
a$.

We now prove (\ref{lem2.1.3}).  Let $p=1-q=\Phi(-z)$. Assume first that $z\ge0$ so that $p\leq 1/2$. $E(T_1|M\geq 1)=t(a,z)$ and for $k\geq 2$,
\begin{eqnarray*}  E(T_k|M\geq k) &=& E(t(a-X_{k-1},z)| M>k-1)\\
&\leq & \mu^{-1} E(a-X_{k-1}| M>k-1)\quad\mbox{since $p\le 1/2$}\\
& \leq& \mu^{-1} g^{(k-1)}(a)
\end{eqnarray*} by (\ref{lem2.1.1}).  Using these two relations $$E(T|M=m) =
\sum_{k=1}^m E(T_k| M=m) = \sum_{k=1}^m E(T_k|M\geq k)\leq t(a,z)+
\mu^{-1} \sum_{k=2}^m g^{(k-1)}(a),$$ since $E(T_k| M=m) = E(T_k|
M\geq k)$ for any $m\geq k$ as discussed above. Thus
\begin{eqnarray*}
  ET &=& E(E(T|M)) \\
   &\leq& t(a,z)+\mu^{-1}\sum_{m\geq 2} q^{m-1}p \sum_{k=2}^m g^{(k-1)}(a)  \\
   &=& t(a,z) + \mu^{-1}\sum_{k\geq 1} g^{(k)}(a)q^k\quad\mbox{by reversing order of summation}\\
   &=& a/\mu + \frac{q\Delta(z)}{\mu\Delta(-z)}\cdot g(a)+\mu^{-1} \sum _{k\geq 2}
   g^{(k)}(a)q^k,
\end{eqnarray*} using the relation between $g$ and $t(\cdot,z)$ in (\ref{lem2.1.4}).

Now assume $z\le0$.  Consequently, $t(\cdot,z)$ is concave, so using Jensen's inequality and (\ref{lem2.1.1}),
\begin{eqnarray*} E(T_k|M\geq k) &=& E[t(a-X_{k-1},z)| M>k-1]\\
& \leq& t(E[a-X_{k-1}| M>k-1],z) \leq t(g^{(k-1)}(a),z)\end{eqnarray*} and, as computed above,
$$ET = E(E(T|M))\le\sum_{m\geq 1}q^{m-1}p \sum_{k=1}^m t(g^{(k-1)}(a),z)=a/\mu + \frac{q\Delta(z)}{\mu\Delta(-z)} \sum _{k\geq 1}g^{(k)}(a)q^{k-1},$$ again using (\ref{lem2.1.4}) for the final step.\qed\bigskip

\noindent {\bf Proof of Lemma \ref{lem2.2}.} Let $z=z(a)$ and $g$ be as in Lemma \ref{lem2.1} with $p=\Phi(-z)$ and $q=1-p$.  A simple computation shows that $g(a)$ has a unique positive fixed point $$x^*=x^*(a)=\frac{\Delta(-z)\phi(z)}{\mu q^2}$$ so that $g(a)\le(a\vee x^*)$, where $\vee$ denotes max.  Then
\begin{eqnarray} EX_M-a&=&\mu(ET-a/\mu)\quad\mbox{Wald's equation}\nonumber\\
&\le&\frac{q\Delta(z)}{\Delta(-z)}\sum_{k=1}^{\infty}g^{(k)}(a)q^{k-1}\quad\mbox{by Lemma \ref{lem2.1}}\nonumber\\
&\le&\frac{q\Delta(z)}{\Delta(-z)}\cdot\begin{cases} \sum_{k=1}^{\infty}g(a)q^{k-1}=g(a)/p\le 2g(a)&\mbox{when $a>x^*$,}\\ \sum_{k=1}^{\infty}x^*q^{k-1}=x^*/p\le 2x^*&\mbox{when $a\le x^*$,}\end{cases}\nonumber\\
&&\quad\quad\mbox{since $g(a)\le(a\vee x^*)$}\nonumber\\
&\le& 2\cdot \frac{q\Delta(z)}{\Delta(-z)}\cdot (g(a)\vee x^*).\label{lem2.2.1}\end{eqnarray} Now \begin{equation}\label{lem2.2.2}\frac{q\Delta(z)}{\Delta(-z)}\cdot g(a)= \frac{\Delta(z)}{2\mu}(\sqrt{4a\mu+z^2}-z)\sim\frac{\abs{z}}{2\mu}\cdot O(\sqrt{a})=O(\abs{z}\sqrt{a})\end{equation} by (\ref{2.8}) and since $\abs{z}=o(\sqrt{a})$. Also,
\begin{eqnarray} \frac{q\Delta(z)}{\Delta(-z)}\cdot x^*&=&\frac{q\Delta(z)}{\Delta(-z)}\cdot \frac{\Delta(-z)\phi(z)}{\mu q^2}=\frac{\Delta(z)\phi(z)}{\mu q}\nonumber\\
&\sim&\frac{|z|\phi(z)}{\mu\Phi(z)}\quad\mbox{by (\ref{2.8}) and since $q=\Phi(z)$}\nonumber\\
&=&O(z^2)\quad\mbox{by (\ref{2.13})}\nonumber\\
&=&O(\abs{z}\sqrt{a})\label{lem2.2.3}\end{eqnarray} since $\abs{z}=o(\sqrt{a})$.  Plugging (\ref{lem2.2.2}) and (\ref{lem2.2.3}) into (\ref{lem2.2.1}) gives the claim.\qed\bigskip

\noindent {\bf Proof of Theorem \ref{thm2.1}.} We proceed by induction on $m$.  For $m=1$, assume $h\in\mathcal{B}_1^o$ and let $\delta_{1,h}^o(a)=(T,M)$.  Note that $$\abs{\zeta(a)}\le \frac{\sqrt{h(a)}}{a^{1/4}}=\frac{\sqrt{h(a)/a}}{a^{1/4}}\cdot\sqrt{a} = o(1)\cdot\sqrt{a}=o(\sqrt{a}),$$ so Lemma \ref{lem2.2} applies.  Then
\begin{eqnarray*} ET-a/\mu&=&\mu\inv(EX_M-a)\quad\mbox{by Wald's equation}\\
&=&\Oh{\abs{\zeta(a)}\sqrt{a}}\quad\mbox{by Lemma (\ref{lem2.2})}\\
&=&\Oh{\frac{\sqrt{h(a)}}{a^{1/4}}\cdot\sqrt{a}}\quad\mbox{by definition of $\zeta$}\\
& =& \Oh{\sqrt{\frac{a^{1/2}}{h(a)}}\cdot h(a)}\\
&=&o(1)\cdot h(a)\quad\mbox{by virtue of $h\in\mathcal{B}_1^o$}\\
&=&o(h(a)),\end{eqnarray*} showing that the theorem holds for $m=1$.

Now assume that $h\in\mathcal{B}_{m+1}^o$ and let $\delta_{m+1,h}^o(a)=(T,M)$.  Let $\xi=\sqrt{\log(a/h(a)^2+1)}$ and $t = t(a,\xi)$, the size of the first stage of $\delta_{m+1,h}^o(a)$.  Note that $$\frac{a}{h(a)^2}=\left(\frac{h_m(a)}{h(a)}\right)^2\cdot \frac{a}{h_m(a)^2}\ge \left(\frac{h_m(a)}{h(a)}\right)^2\cdot \frac{a}{h_1(a)^2}=\left(\frac{h_m(a)}{h(a)}\right)^2\cdot 1\To\infty$$ so that $\xi\To\infty$, which is important in the definition of the first stage of $\delta_{m+1,h}^o$.  Using the definition of $f$, as $x\To\infty$ \begin{eqnarray*} h_m(f(x))&=& O((x\log x)^{(1/2)^{m+1}}(\log (x\log x))^{1/2-(1/2)^m})\\
&=&O(x^{(1/2)^{m+1}}(\log x)^{1/2-(1/2)^{m+1}})=O(h_{m+1}(x))=o(h(x)),\end{eqnarray*} or equivalently, $h_m=o(h\circ f\inv)$.  A similar argument shows $h\circ f\inv=o(h_{m-1})$, which shows that \begin{equation}\label{thm2.1.4}h\circ f\inv\in\mathcal{B}_m^o.\end{equation} For $y>0$ let $\delta_{m,h\circ f\inv}^o(y)=(T',M')$ and 
\begin{eqnarray*} \psi(y)&=&EM'\\
\varphi(y)&=&EX_{M'}-y.\end{eqnarray*}  Suppose $\eps>0$.  By the induction hypothesis there is a constant $y_o$ such that 
\begin{eqnarray} \abs{\psi(y)-m}&\le&\eps\label{thm2.1.10}\\
\varphi(y)&\le&\eps\cdot h(f\inv(y))\label{thm2.1.6}\end{eqnarray} for all $y\ge y_o$, where (\ref{thm2.1.6}) uses (\ref{thm2.1.1}) and Wald's equation.
\begin{eqnarray} EM-1&=&E(M-1;0<a-X_1<y_o)+E(M-1;a-X_1\ge y_o)\nonumber\\
&=&E(\psi(a-X_1);0<a-X_1<y_o)+E(\psi(a-X_1);X_1\ge y_o).\label{thm2.1.12}\end{eqnarray}  Letting $A=\sup_{0<y<y_o}\psi(y)<\infty$,
\begin{eqnarray} E(\psi(a-X_1);0<a-X_1<y_o)&\le&A\cdot P(0<a-X_1<y_o)\nonumber\\
&\le&A\cdot P(X_1>a-y_o)\nonumber\\
&=&A\cdot P\left(\frac{X_1-\mu t}{\sqrt{t}}>\frac{a-\mu t}{\sqrt{t}}-\frac{y_o}{\sqrt{t}}\right)\nonumber\\
&=&A\cdot P\left(\frac{X_1-\mu t}{\sqrt{t}}>\xi+o(1)\right)\nonumber\\
&=&A\cdot o(1)=o(1),\label{thm2.1.11}\end{eqnarray} by, say, Chebyshev's inequality since $\xi\To\infty$.  Plugging this into (\ref{thm2.1.12}),
\begin{eqnarray*} \abs{EM-(m+1)}&\le&o(1)+\abs{E(\psi(a-X_1)-m;a-X_1\ge y_o)}+m\cdot P(a-X_1\le y_o)\\
&\le& o(1)+E(\abs{\psi(a-X_1)-m};a-X_1\ge y_o)+m\cdot o(1)\\
&\le&o(1)+\eps\cdot P(a-X_1\ge y_o)+o(1)\quad\mbox{by (\ref{thm2.1.10})}\\
&\le& \eps+\eps+\eps=3\eps\end{eqnarray*} for sufficiently large $a$.  This shows that $EM\To m+1$.

To show that (\ref{thm2.1.2}) holds for $m+1$, first write 
\begin{eqnarray}\lefteqn{ \mu E(T-a/\mu)=EX_M-a}\nonumber\\
&&=E(X_M-a;M=1)+E(X_M-a;0<a-X_1<y_o)\nonumber\\
&&\quad+E(X_M-a;y_o\le a-X_1\le (1/2)f(a))+E(X_M-a;a-X_1>(1/2) f(a))\nonumber\\
&&\equiv I+II+III+IV.\label{thm2.1.17}\end{eqnarray} 

First,
\begin{eqnarray} I=E(X_1-a;X_1\ge a)&=& \sqrt{t}\cdot \Delta(\xi)\quad\mbox{by (\ref{2.11})}\nonumber\\
&\sim& \sqrt{t}\cdot \frac{\phi(\xi)}{\xi^2}\quad\mbox{by (\ref{2.8})}\nonumber\\
&=&O(\sqrt{a})\cdot\frac{O(h(a)/\sqrt{a})}{\xi^2}\nonumber\\
&&\mbox{(since $t\sim a/\mu$ and $\phi(\xi)\propto h/\sqrt{a}$)}\nonumber\\
&=&o(h(a)).\label{thm2.1.13}\end{eqnarray} Next we have \begin{equation}\label{thm2.1.14} II=E(\varphi(a-X_1);0<a-X_1<y_o)=o(1)\end{equation} by the same argument leading to (\ref{thm2.1.11}).

Before considering $III$, note that we may assume without loss of generality that $h$ is non-decreasing.  Otherwise, we could replace $h$ by $\underline{h}(x)\equiv\inf_{y\ge x}h(y)$ throughout the proof, since $\underline{h}$ is non-decreasing and bounded above by $h$.  Since $f\inv$ is also non-decreasing, $h\circ f\inv$ is thus non-decreasing. Now 
\begin{eqnarray} III&=&E(\varphi(a-X_1);y_o\le a-X_1\le(1/2) f(a))\nonumber\\
&\le&E(\eps\cdot h(f\inv(a-X_1));y_o\le a-X_1\le(1/2) f(a))\quad\mbox{by (\ref{thm2.1.6})}\nonumber\\
&\le&\eps\cdot h(f\inv((1/2) f(a)))\quad\mbox{$h\circ f\inv$ non-decreasing}\nonumber\\
&\le&\eps\cdot h(f\inv(f(a)))\quad\mbox{$f\inv$ non-decreasing}\nonumber\\
&=&\eps\cdot h(a).\label{thm2.1.15}\end{eqnarray}

We know that $h(f\inv(x))=o(h_{m-1}(x))$ and hence $h(f\inv(x))=o(x)$ since $m=1$ gives the largest case, asymptotically.  Thus we assume that $a$ is large enough so that \begin{equation}\label{thm2.1.7}h(f\inv(y))\le y\quad\mbox{for all $y\ge (1/2)f(a)$.}\end{equation} Then 
\begin{eqnarray} IV&=&E(\varphi(a-X_1);a-X_1>(1/2)f(a))\nonumber\\
&\le&E(\eps\cdot h(f\inv(a-X_1));a-X_1>(1/2)f(a))\quad\mbox{by (\ref{thm2.1.6})}\nonumber\\
&\le&\eps E(a-X_1;a-X_1>(1/2)f(a))\quad\mbox{by (\ref{thm2.1.7})}\nonumber\\
&=&\eps\left[\sqrt{t}\cdot \Delta\left(\frac{f(a)}{2\sqrt{t}}-\xi\right)+(f(a)/2)\cdot\Phi\left(\xi-\frac{f(a)}{2\sqrt{t}}\right)\right],\label{thm2.1.9}\end{eqnarray} using (\ref{2.12}) for this last step.  Now 
\begin{eqnarray*} \frac{f(a)}{2\sqrt{t}}-\xi&=&-\xi+\frac{f(a)}{2\sqrt{a/\mu}}\cdot(1+o(1))\quad\mbox{since $t\sim a/\mu$}\\
&=&-\sqrt{\log(a/h(a)^2)}+o(1)+3\sqrt{\log(a+1)}\cdot (1+o(1))\\
&\ge&-\sqrt{\log(a/h_{m+1}(a)^2)}+2\sqrt{\log a}\quad\mbox{since $h\in\mathcal{B}_m^o$}\\
&\ge&-\sqrt{\log(a/a^{(1/2)^m})}+2\sqrt{\log a}\\
&=&-\sqrt{(1-2^{-m})\log a}+2\sqrt{\log a}\ge \sqrt{\log a}.\end{eqnarray*} Plugging this back into (\ref{thm2.1.9}) gives 
\begin{eqnarray} IV&\le&\eps\left[\sqrt{t}\cdot\Delta(\sqrt{\log a})+(f(a)/2)\cdot\Phi(-\sqrt{\log a})\right]\nonumber\\
&=&O(\sqrt{a})\cdot\frac{\phi(\sqrt{\log a})}{\log a}+O(\sqrt{a\log a})\cdot \frac{\phi(\sqrt{\log a})}{\sqrt{\log a}}\quad\mbox{by (\ref{2.8}) and Mills' ratio}\nonumber\\
&=&O(\sqrt{a})\cdot\phi(\sqrt{\log a})=O(1)=o(h(a)).\label{thm2.1.16}\end{eqnarray} Plugging (\ref{thm2.1.13}), (\ref{thm2.1.14}), (\ref{thm2.1.15}), and (\ref{thm2.1.16}) into (\ref{thm2.1.17}) gives $$\mu ET-a/\mu\le o(h(a))+o(1)+\eps h(a)+o(h(a))\le 4\eps h(a)$$ for sufficiently large $a$.  This shows that $ET-a/\mu=o(h(a))$, finishing the proof.\qed\bigskip

\noindent {\bf Proof of Lemma \ref{lem2.3}.} \begin{eqnarray*} \log(\kappa_{m+1}/\kappa_m)&=&-(1/2)^{m+1}\log\mu+\sum_{i=1}^m(1/2)^{i+1}\log[(1/2)^{m-i}-(1/2)^m]\\
&&-\sum_{i=1}^{m-1}(1/2)^{i+1}\log[(1/2)^{m-1-i}-(1/2)^{m-1}]\\
&=&-(1/2)^{m+1}\log\mu +(1/2)^{m+1}\log[1-2^{-m}]-\sum_{i=1}^{m-1}(1/2)^{i+1}\log 2\\
&=&(1/2)^{m+1}\log(1-2^{-m})/\mu+(1/2-(1/2)^m)\log (1/2).\end{eqnarray*} On the other hand, letting $y=\sqrt{(1-2^{-m})/\mu\cdot x\log x}$,
\begin{eqnarray*} \log(h_m(y)/h_{m+1}(x))&=&(1/2)^m\log y+(1/2-(1/2)^m)\log\log y\\
&&-(1/2)^{m+1}\log x-(1/2-(1/2)^{m+1})\log\log x\\
&=&(1/2)^{m+1}[\log (1-(1/2)^m)/\mu+\log x+\log\log x]\\
&&+(1/2-(1/2)^m)[\log (1/2)+\log\log x+o(1)]\\
&&-(1/2)^{m+1}\log x-(1/2-(1/2)^{m+1})\log\log x\\
&=&(1/2)^{m+1}\log(1-2^{-m})/\mu+(1/2-(1/2)^m)\log (1/2)+o(1)\\
&=&\log(\kappa_{m+1}/\kappa_m)+o(1)\end{eqnarray*} so that $h_m(y)/h_{m+1}(x)\To \kappa_{m+1}/\kappa_m$.\qed\bigskip

\noindent {\bf Proof of Theorem \ref{thm2.2}.} We proceed by induction on $m$.  Let $\delta_{1,z}^+(a)=(T,M)$.  Using Wald's equation, $$\mu E(T-a/\mu)=E(X_M-a;M=1)+E(X_M-a;M>1).$$  Let $t=t(a,z)$, the size of the first stage of $\delta_{1,z}^+(a)$.  Using (\ref{2.8}), 
\begin{eqnarray} E(X_M-a;M=1)&=&E(X_1-a;X_1\ge a)\nonumber\\
&=&\sqrt{t}\cdot\Delta(z)\quad\mbox{by (\ref{2.11})}\nonumber\\
&\sim&\sqrt{a/\mu}\cdot \Delta(z)\label{thm2.2.3} \end{eqnarray} since $t\sim a/\mu$.  For $y>0$ let \begin{eqnarray*}\dot{\delta}_{\nu(y)}(y)&=&(T',M'),\\
\varphi(y)&=&EX_{M'}-y,\\
\psi(y)&=&EM',\end{eqnarray*} where $\nu(y)=-\sqrt{\log(y+1)}$ as in the definition of $\delta_{1,z}^+$.  Suppose $\eps>0$.  By Lemma \ref{lem2.2} and since $EM'=1/\Phi(|\nu(y)|)\To1$ as $y\To\infty$ there are constants $y_o,C$ such that for all $y\ge y_o$,
\begin{eqnarray} \varphi(y)&\le& C\sqrt{y\log(y+1)}\label{thm2.2.4}\\
\abs{\psi(y)-1}&\le&\eps.\label{thm2.2.12}\end{eqnarray}  Let $Y=a-X_1$.  Note that 
\begin{equation} P(Y\ge y_o)=\Phi(z-y_o/\sqrt{t})\To \Phi(z)\label{thm2.2.5}\end{equation} by continuity, and  \begin{equation}\label{thm2.2.6}P(Y<0)= 1-\Phi(z)\end{equation} so that $P(0<Y<y_o)\To 0$.

Now
\begin{eqnarray*} E(X_M-a;M>1)&=&E(\varphi(Y);Y>0)\\
&=&E(\varphi(Y);0<Y<y_o)+E(\varphi(Y);Y\ge y_o)\\
&=&o(1)+E(C\sqrt{Y\log(Y+1)};Y\ge y_o)\quad\mbox{by \ref{thm2.2.6}}\\
&\le&o(1)+C\sqrt{E(Y|Y\ge y_o)\log [E(Y|Y\ge y_o)+1]},\end{eqnarray*} where this last uses concavity of $y\mapsto y\log(y+1)$ with Jensen's inequality.  Now  
\begin{eqnarray*} \lefteqn{E(Y|Y\ge y_o)=P(Y\ge y_o)\inv E(Y;Y\ge y_o)}\\
&=&[\Phi(z)+o(1)]\inv\cdot [\sqrt{t}\cdot\Delta(-z+y_o/\sqrt{t})+y_o\cdot\Phi(z-y_o/\sqrt{t})]\quad\mbox{using (\ref{2.12})}\\
&=&O(1)[O(\sqrt{a})\cdot O(1)+y_o\cdot O(1)]=O(\sqrt{a}).\end{eqnarray*}  Thus \begin{eqnarray*} E(X_M-a;M>1)&\le&o(1)+C\sqrt{O(\sqrt{a})\log O(\sqrt{a})}\\
&=&o(\sqrt{a}),\end{eqnarray*} and combining this with (\ref{thm2.2.3}) gives $$ET-a/\mu\le \Delta(z)\mu^{-3/2}\sqrt{a}+o(\sqrt{a})=(1+o(1))\Delta(z)\kappa_1 h_1(a),$$ establishing (\ref{thm2.2.1}) for $m=1$.
\begin{eqnarray*}\lefteqn{\abs{EM-(1+\Phi(z))}=\abs{1+E(\psi(Y);Y>0)-(1+\Phi(z))}}\\
&\le&\abs{E(\psi(Y);0<Y<y_o)}+E(\abs{\psi(Y)-1};Y\ge y_o)+\abs{P(Y\ge y_o)-\Phi(z)}\\
&\le&o(1)+\eps+o(1)\quad\mbox{by (\ref{thm2.2.12}) and (\ref{thm2.2.5})}\\
&\le&3\eps\end{eqnarray*} for sufficiently large $a$, proving that $EM\To 1+\Phi(z)$.

Next let $\delta_{m+1,z}^+(a)=(T,M)$ and $\delta_{m,z}^+(y)=(T',M')$ for $y>0$.  Let 
\begin{eqnarray*} \varphi(y)&=&EX_{M'}-y,\\
\psi(y)&=&EM',\end{eqnarray*} and suppose $\eps>0$.  By the induction hypothesis there is a constant $y_1$ such that, for all  $y\ge y_1$,
\begin{eqnarray} \varphi(y)&\le&(1+\eps)\mu\Delta(z)\kappa_mh_m(y)\label{thm2.2.7}\\
\abs{\psi(y)-(m+\Phi(z))}&\le&\eps\label{thm2.2.8}.\end{eqnarray} Let $t=t(a,\sqrt{(1-2^{-m})\log(a+1)})$, the size of the first stage of $\delta_{m+1,z}^+(a)$.  Again letting $Y=a-X_1$,
\begin{eqnarray} \mu E(T-a/\mu)&=&E(X_M-a;M=1)+E(X_M-a;M>1)\nonumber\\
&=&E(X_1-a;X_1\ge a)+E(\varphi(Y);0<Y<y_1)+E(\varphi(Y);Y\ge y_1).\label{thm2.2.9}\end{eqnarray} Now
\begin{eqnarray} E(X_1-a;X_1\ge a)&=&\sqrt{t}\cdot\Delta(\sqrt{(1-2^{-m})\log(a+1)})\quad\mbox{by (\ref{2.11})}\nonumber\\
&\sim&\sqrt{a/\mu}\cdot \frac{\phi(\sqrt{(1-2^{-m})\log(a+1)})}{(1-2^{-m})\log(a+1)}\nonumber\quad\mbox{by (\ref{2.8})}\\
&=&\sqrt{a/\mu}\cdot\frac{O(a^{-1/2+(1/2)^{m+1}})}{(1-2^{-m})\log(a+1)}\nonumber\\
&=&o(a^{(1/2)^{m+1}})=o(h_{m+1}(a)).\label{thm2.2.10}\end{eqnarray} $E(\varphi(Y);0<Y<y_1)=o(1)$ by a now routine argument and using (\ref{thm2.2.7}),
\begin{eqnarray*} \frac{E(\varphi(Y);Y\ge y_1)}{(1+\eps)\mu\Delta(z)\kappa_m}&\le& E(h_m(Y);Y\ge y_1)\\
&\le&h_m(E(Y;Y\ge y_1)P(Y\ge y_1)\inv)\end{eqnarray*} by Jensen's inequality since $h_m$ is concave.  Now $$E(Y;Y\ge y_1)\sim\sqrt{t}\cdot \Delta(-\sqrt{(1-2^{-m})\log(a+1)})\sim\sqrt{a/\mu}\cdot \sqrt{(1-2^{-m})\log(a+1)}$$ by a routine application of (\ref{2.12}).  Using that $P(Y\ge y_1)\To1$, for sufficiently large $a$ we have
\begin{eqnarray} E(\varphi(Y);Y\ge y_1)&\le&(1+\eps)^2\mu\Delta(z)\kappa_mh_m(\sqrt{a/\mu})\nonumber\\
&\le&(1+\eps)^3\mu\Delta(z)\kappa_{m+1}h_{m+1}(a),\label{thm2.2.11}\end{eqnarray} using Lemma \ref{lem2.3}.  Plugging (\ref{thm2.2.10}) and (\ref{thm2.2.11}) into (\ref{thm2.2.9}) gives
\begin{eqnarray*} ET-a/\mu&\le&o(h_{m+1}(a))+o(1)+(1+\eps)^3\Delta(z)\kappa_{m+1}h_{m+1}(a)\\
&\le&[\eps+\eps+(1+\eps)^3]\Delta(z)\kappa_{m+1}h_{m+1}(a)\\
&\le&[1+9\eps]\Delta(z)\kappa_{m+1}h_{m+1}(a)\end{eqnarray*} for sufficiently large $a$, which shows that (\ref{thm2.2.1}) holds for $m+1$.  As for the number of stages used,
\begin{eqnarray*} \lefteqn{\abs{EM-(m+1+\Phi(z))}=\abs{1+E(\psi(Y);Y>0)-(m+1+\Phi(z))}}\\
&\le&\abs{E(\psi(Y);0<Y< y_1)}+E(\abs{\psi(Y)-(m+\Phi(z))};Y\ge y_1)+(m+\Phi(z))P(Y<y_1)\\
&\le&o(1)+\eps\cdot P(Y\ge y_1)+(m+\Phi(z))\cdot o(1)\quad\mbox{by (\ref{thm2.2.8}) and since $P(Y<y_1)\To 0$}\\
&\le&o(1)+\eps+o(1)\le 3\eps\end{eqnarray*} for sufficiently large $a$.  This shows that $EM\To m+1+\Phi(z)$ and completes the proof.\qed\bigskip

\noindent {\bf Proof of Lemma \ref{lem2.5}.}  Let $G_k = (1/\mu)^{1-2^{-k}}
F_h^{(k)}(a)$.  Suppose $\varepsilon>0$.  Let $$V_k = \{a-X_k\geq
(1-\varepsilon)G_k\}.$$  The $k=0$ case is trivial since
(\ref{lem2.5.0}) is equivalent to $P(a\ge a)\To 1$.  Fix $1\leq k<m$ and assume
that $P(V_{k-1})\To 1$.  Let $\delta(a)=(T,M)$ and
$$\zeta_k = \frac{a-X_{k-1}-\mu T_k}{\sqrt{T_k}}.$$ Note that
\begin{eqnarray*} h(a)^2 &=& o(h_{m-1}(a)^2)\quad\mbox{since $h\in\mathcal{B}_m$}\\
& =&o(F_h^{(m-2)}(a))\quad\mbox{by Lemma \ref{lem2.4}}\\
&=&o(G_{m-2})\\
&=&o(G_{k-1})\end{eqnarray*} since $k-1\leq m-2$.  Thus $G_{k-1}/h^2\To\infty$ and so does $\log(G_{k-1}/h^2)$.  With this, we claim
\begin{equation} \label{lem2.5a} P(\zeta_k \geq
\sqrt{\log(G_{k-1}/h^2(a))}-1 |V_{k-1})\To 1.
\end{equation}Let $\xi = \sqrt{\log(G_{k-1}/h^2(a))}-1$ and $U = \{\zeta_{k} < \xi\}$.  If (\ref{lem2.5a})
were to fail there would be a constant $\eta>0$ and a sequence of $a$'s approaching $\infty$ on which $P(U|V_{k-1})>\eta$.  Then \begin{eqnarray}
\mu R(\delta(a)) & \geq & \mu E(T-a/\mu) =E(X_M-a) \geq E[(X_M-a)1\{M=k\}; U\cap V_{k-1}]\nonumber \\
&=& E[\Delta(\zeta_{k}) \sqrt{t(a-X_{k-1},\zeta_k)}; U\cap V_{k-1}].\label{lem2.5b}
\end{eqnarray}The function inside the expectation in (\ref{lem2.5b}) is decreasing
in both $\zeta_{k}$ and $X_{k-1}$, hence
\begin{equation}
\label{lem2.5e} \mu R(\delta(a))\geq \Delta(\xi)
\sqrt{t((1-\varepsilon)G_{k-1},\xi)}\cdot P(U\cap V_{k-1}).
\end{equation}  By assumption, $P(U|V_{k-1})\geq \eta$ and $P(V_{k-1})\To 1$, so \begin{equation}
\label{lem2.5c}
P(U\cap V_{k-1})\geq \eta/2,
\end{equation} say, for large enough $a$.  Also
\begin{eqnarray} \lefteqn{\Delta(\xi) \sqrt{t((1-\varepsilon)G_{k-1},\xi)} \sim  \frac{\phi(\xi)}{\xi^2} \sqrt{(1-\varepsilon)G_{k-1}/\mu}}\nonumber\\
 &\geq& \varepsilon'  h(a)\frac{\exp(\sqrt{\log(G_{k-1}/h^2(a))}-1/2)}{(\sqrt{\log(G_{k-1}/h^2(a))}-1)^2}\quad\mbox{some $\varepsilon'>0$}\nonumber\\
 &=& h(a)/o(1). \label{lem2.5d}\end{eqnarray} Plugging (\ref{lem2.5c}) and (\ref{lem2.5d}) into (\ref{lem2.5e}) gives $h(a)=o(R(\delta(a)))$, which contradicts our assumption that $R(\delta(a))=O(h(a))$.  Hence, (\ref{lem2.5a}) must hold.  Then
\begin{equation}\begin{split}
P(& V_{k}| U'\cap V_{k-1}) =P(a-X_{k}\geq (1-\varepsilon)G_{k}|U'\cap V_{k-1}) \\
&= P\left(\left.\frac{(X_{k}-X_{k-1}) - \mu T_{k}}{\sqrt{T_{k}}}\leq \frac{a-X_{k-1}-  (1-\varepsilon)G_{k}-\mu T_{k}}{\sqrt{T_{k}}}\right| U'\cap V_{k-1}\right).\label{lem2.5f} \end{split}\end{equation}On $V_{k-1}$,
\begin{eqnarray*}
 \frac{a-X_{k-1}-  (1-\varepsilon)G_{k}-\mu T_{k}}{\sqrt{T_{k}}} & = &  \zeta_{k} -\frac{2\mu (1-\varepsilon)G_{k}}{\sqrt{4\mu (a-X_{k-1})+ \zeta_{k}^2}-\zeta_{k}} \\
 & \geq & \zeta_{k} -\frac{2\mu (1-\varepsilon)G_{k}}{\sqrt{4\mu (1-\varepsilon)G_{k-1}+ \zeta_{k}^2}-\zeta_{k}},
\end{eqnarray*}which is increasing in $\zeta_{k}$.  Hence, on $U'$,
\begin{equation*}\begin{split}
\zeta_{k} -&\frac{2\mu (1-\varepsilon)G_{k}}{\sqrt{4\mu (1-\varepsilon)G_{k-1}+ \zeta_{k}^2}-\zeta_{k}} \geq \xi -\frac{2\mu (1-\varepsilon)G_{k}}{\sqrt{4\mu (1-\varepsilon)G_{k-1}+\xi^2}-\xi}\\
&=\xi-\frac{2\mu (1-\varepsilon)G_{k}}{\sqrt{4\mu (1-\varepsilon)G_{k-1}}}(1+o(1))\\
&= \xi-\sqrt{(1-\varepsilon)\log(F_h^{(k-1)}(a)/h^2(a))}(1+o(1)) \\
&\sim (1-\sqrt{1-\varepsilon}) \sqrt{\log(G_{k-1}/h^2(a))}\equiv\gamma\To\infty.\end{split}\end{equation*}  Substituting this back into (\ref{lem2.5f}) gives $$P(V_{k}|U'\cap V_{k-1})\ge1-(\gamma/2)^{-2}\To 1$$ by Chebyshev's inequality.  Thus $P(V_{k})\geq P(V_{k}|U'\cap V_{k-1}) P(U'\cap V_{k-1})\To 1$ since  $P(U'\cap V_{k-1})\To 1$ by the induction hypothesis and (\ref{lem2.5a}), finishing the induction and proving the lemma.\qed\bigskip

\noindent {\bf Proof of Lemma \ref{lem2.4}.} Let $F^k$ denote $F_{h}^{(k)}(x)$.  First we prove (\ref{lem2.4.2}) by induction on $k$.  For $k=1$, $$\sqrt{F^0}=\sqrt{x}=1\cdot \sqrt{x}=C_1^m\cdot h_1(x).$$  Now assume $2\le k+1\le m$, $\sqrt{F^{k-1}}\sim C_k^mh_k(x)$, and let $Q=\lim h/h_m\in(0,\infty)$.  Observe that
\begin{eqnarray} \log\left(\frac{F^{k-1}}{h(x)^2}\right)&\sim& \log\left(\frac{(C_k^mh_k(x))^2}{(Qh_m(x))^2}\right)\nonumber\\
&\sim& \log\left(\frac{h_k(x)^2}{h_m(x)^2}\right)\nonumber\\
&\sim& \log\left(\frac{x^{(1/2)^{k-1}}(\log x)^{1-(1/2)^{k-1}}}{x^{(1/2)^{m-1}}(\log x)^{1-(1/2)^{m-1}}}\right)\nonumber\\
&\sim&[(1/2)^{k-1}-(1/2)^{m-1}]\log x,\label{lem2.4.4}\end{eqnarray} so
\begin{eqnarray} \sqrt{F^k}&=&\set{F^{k-1}\log (F^{k-1}/h(x)^2)}^{1/4}\nonumber\\
&\sim&\set{(C_k^mh_k(x))^2[(1/2)^{k-1}-(1/2)^{m-1}]\log x}^{1/4}\nonumber\\
&=&\sqrt{C_k^m}\cdot x^{(1/2)^{k+1}}(\log x)^{1/4-(1/2)^{k+1}}[(1/2)^{k-1}-(1/2)^{m-1}]^{1/4}(\log x)^{1/4}\nonumber\\
&=&\sqrt{C_k^m}\cdot [(1/2)^{k-1}-(1/2)^{m-1}]^{1/4}h_{k+1}(x)\nonumber\\
&=&C_{k+1}^mh_{k+1}(x),\label{lem2.4.5}\end{eqnarray} by (\ref{lem2.4.6}).

Next we prove (\ref{lem2.4.1}) by induction on $k$.  The $k=1$ case is again easy since $$C_1^{m-1}=\frac{\sqrt{F^0}}{h_1}=C_1^m=1$$ for any $m\ge 2$.  Now assume $2\le k+1<m$ and that  (\ref{lem2.4.1}) holds for $k$.  Then, since $h_m\ll h\ll h_{m-1}$,
$$\log\left(\frac{F^{k-1}}{h(x)^2}\right)\lesssim\log\left(\frac{(C_k^mh_k(x))^2}{h_m(x)^2}\right)\sim [(1/2)^{k-1}-(1/2)^{m-1}]\log x,$$ by the same argument leading to (\ref{lem2.4.4}).  Then, by repeating the argument leading to (\ref{lem2.4.5}) with $\lesssim$ in place of $\sim$, $$\sqrt{F^k}\lesssim \sqrt{C_k^m}\cdot [(1/2)^{k-1}-(1/2)^{m-1}]h_{k+1}(x)=C_{k+1}^mh_{k+1}(x),$$ by (\ref{lem2.4.6}).  The other bound is similar: $$\log\left(\frac{F^{k-1}}{h(x)^2}\right)\gtrsim\log\left(\frac{(C_k^{m-1}h_k(x))^2}{h_{m-1}(x)^2}\right)\sim [(1/2)^{k-1}-(1/2)^{m-2}]\log x,$$ and so $$\sqrt{F^k}\gtrsim \sqrt{C_k^{m-1}}\cdot [(1/2)^{k-1}-(1/2)^{m-2}]h_{k+1}(x)=C_{k+1}^{m-1}h_{k+1}(x),$$ by replacing $m$ by $m-1$ in (\ref{lem2.4.5}) and (\ref{lem2.4.6}).\qed\bigskip

\section*{ACKNOWLEDGEMENTS}The author would like to thank Gary Lorden
for suggesting this problem and for his enthusiastic support and
guidance of the author's thesis research, of which this work is a
part.  This work is supported by NSF DMS-0403105.

\section*{REFERENCES}
\begin{list}{}{\setlength{\labelwidth}{0cm}\setlength{\labelsep}{0cm}\setlength{\leftmargin}{1cm}\setlength{\itemindent}{-1cm}}
\item{}Bartroff, J. (2004). Asymptotically Optimal Multistage Hypothesis Tests, Ph.D.\ diss., Caltech.
\item{}Bartroff, J. (2005). Asymptotically optimal multistage tests of simple hypotheses, submitted.
\item{}Cressie, N.\ and Morgan P.B. (1993). The VPRT: A sequential
  testing procedure dominating the SPRT, \textit{Econometric Theory}
  9: 431--450.
\item{}Feller, W. (1971). \textit{An Introduction to Probability
    Theory and its Applications}, vol. 2, second edition, New York: Wiley.
\item{}Lorden, G. (1970). On excess over the boundary,
  \textit{Annals of Mathematical Statistics} 41: 520-527.
\item{}Morgan, P. B. and Cressie, N. (1997). A comparison of the
  cost-efficiencies of the sequential, group-sequential, and
  variable-sample-size-sequential probability ratio tests,
  \textit{Scandinavian Journal of Statistics} 24: 181-200.
\item{}Schmitz, N. (1993). \textit{Optimal sequentially planned
    decision procedures}, Lecture Notes in Statistics, vol.\ 79, New
    York: Springer-Verlag.
\item{}Siegmund, D. (1985). \textit{Sequential Analysis}, New York:
  Springer-Verlag. 
\end{list}


\end{document}